\documentclass[10pt,a4paper]{article}

\usepackage{amsmath}
\usepackage{amssymb}
\usepackage{theorem}
\theoremheaderfont{\scshape}
\newtheorem{thm}{Theorem}[section]
\newtheorem{cor}[thm]{Corollary}
\newtheorem{lem}[thm]{Lemma}

\theorembodyfont{\normalfont}

\def\bb#1{\mathbb{#1}}
\def\cal#1{\mathcal{#1}}
\def\ds{\displaystyle}
\def\eps{\varepsilon}
\def\g#1{\mathbf{#1}}

\def\limf#1#2{\ensuremath{\underset{#1\rightarrow #2}{\lim} }}
\def\abs#1{\left\vert #1\right\vert}
\def\N#1{\|#1\|}
\newcommand{\Lim}[2]{\lim\limits_{#1\rightarrow #2}}

\title{Concentration inequalities for disordered models }
\author{Fr\'ed\'erique Watbled\\
LMBA, UMR 6205, Univ.  Bretagne Sud, \\
Campus de Tohannic, BP 573 \\
 56017 Vannes, France.\\
 frederique.watbled@univ-ubs.fr}
 
\date{}

\begin{document}

\maketitle

\begin{abstract}
We use a generalization of Hoeffding's inequality to show concentration results for the free energy of disordered pinning models, assuming only that the disorder has a finite exponential moment. We also prove some concentration inequalities for directed polymers in random environment, which we use to establish a large deviations result for the end position of the polymer under the polymer measure.
\end{abstract}


%secondary 60K37, 60G50, 60F10, 82D30

%%%%%%% section 
\section{Introduction}
\setcounter{equation}{0}
Let $S=(S_n)_{n\in\bb N}$ be a symmetric random walk defined on a probability space
$(\Sigma,\cal E,\textbf{P})$ with values in $\bb Z^d$, starting at 0.
Let $I$ be a denumerable set, and
let $\eta=(\eta_k)_{k\in I}$ be a sequence of i.i.d. random variables
defined on another probability space $(\Omega,\cal F, \bb P)$.
The expectation of a function $f$ with respect to the probability measure 
$\textbf{P}$ (respectively $\bb P$) will be denoted by $\textbf{E}f=\int_\Sigma f d\textbf{P}$
(respectively $\bb E f=\int_\Omega f d\bb P$).
We assume that there exists $\beta>0$ such that $\bb E e^{\beta|\eta_i|}<\infty$ for $i\in I$.
For every fixed realization of $\eta$,
we consider functionals of the form 
\begin{equation}\label{partition}
Y_n^\eta=\textbf{E}\exp(H_n^\eta(S))f_n(S_n),
\end{equation}
where $H_n^\eta$ is a function of the first $n$ terms of $S$, and $f_n$ is a bounded positive function defined on $\bb R^d$ with $\g P(f_n(S_n)>0)\neq 0$.
When $$I=\bb N,\ H_n^\eta(S)=\sum_{k=1}^n(\beta \eta_k+h)\textbf{1}_{\{S_k=0\}},
\textrm { and } f_n(x)=\textbf{1}_{\{x=0\}},$$ 
$Y_n^\eta$ is the partition function of the disordered pinning model studied in \cite{Giacomin2007}$, \cite{Giacomin2010}$.
When $$I=\bb N\times\bb Z^d,\ H_n^\eta(S)=\beta\sum_{j=1}^n \eta_{(j,S_j)}, \textrm{ and } f_n=1,$$ 
$Y_n^\eta$
is the partition function of a directed polymer in a random environment. 
In this paper we use the exponential inequality for martingales that we obtained in \cite{LiuWatbled}
to prove some concentration properties of the free energy for both models.

Let us present the two models in more detail. We refer to \cite{Giacomin2007}, \cite{Giacomin2010}, for a survey on
 disordered pinning models. 
Let $S=(S_n)_{n\in\bb N}$ be a random walk on
$\bb Z^d$ starting at $0$, defined on a probability space
$(\Sigma,\cal E,\textbf{P})$. We suppose that the increment variables 
$(S_n-S_{n-1})_{n\geq 1}$ are independent, symmetric, and have finite variance.
For each $\beta\geq 0$, $h\in\bb R$ and 
$\eta=(\eta_k)_{k\in\bb N}$ a sequence of real numbers, we define the Gibbs measure $\textbf{P}_n$ on $\Sigma$
by giving its density:
\begin{equation}
\frac{d\textbf{P}_n}{d\textbf{P}}=\ds\frac{1}{Z_n}
\exp(\sum_{k=1}^n(\beta \eta_k+h)\textbf{1}_{\{S_k=0\}})\textbf{1}_{\{S_n=0\}},
\end{equation}
where $Z_n$ is the partition function
\begin{equation}
Z_n=\textbf{E}\exp(\sum_{k=1}^n(\beta \eta_k+h)\textbf{1}_{\{S_k=0\}})\textbf{1}_{\{S_n=0\}}.
\end{equation}
The set $\tau=\{n\geq 0: S_n=0\}$ of visits to 0 is then a renewal process.
It means that if $\tau=\{\tau_0,\tau_1,\cdots\}$ with $\tau_0=0$ and $\tau_i>\tau_{i-1}$, then
$(\tau_{i}-\tau_{i-1})_{i\geq 1}$ is an i.i.d. sequence of positive random variables.
Let $\delta_n=\textbf{1}_{\{n\in\tau\}}$.
Then we can write
\begin{equation}\label{renewal}
\frac{d\textbf{P}_n}{d\textbf{P}}=\ds\frac{1}{Z_n}\exp(\sum_{k=1}^n(\beta \eta_k+h)\delta_k)\delta_n,
\end{equation}
and 
\begin{equation}
Z_n=\textbf{E}\exp(\sum_{k=1}^n(\beta \eta_k+h)\delta_k)\delta_n,
\end{equation}
without mentioning trajectories $S$ anymore.
In the following we consider a renewal process
$\tau=(\tau_j)_{j\in\bb N}$ on a probability space $(\Sigma,\cal E,\textbf{P})$ with values in $\bb N$,
starting with $\tau_0=0$, and
a realization $\eta$ of an i.i.d. sequence of random variables on a 
probability space $(\Omega,\cal F, \bb P)$. For simplicity we write $\bb E f(\eta)=\bb E f(\eta_0)$ for any $f$ such that $f\circ \eta_0$ is integrable. We fix $\beta >0$ and assume that $\bb E e^{\beta|\eta|}<\infty$.  We denote by $\lambda(\beta)=\ln \bb E e^{\beta\eta}$ the logarithmic moment generating function of
$\eta$. 
Thanks to subadditivity properties it is easily shown that $\frac{1}{n}\bb E\ln Z_n$
converges to a number $F(\beta,h)$ which is called the free energy of the model, and it is also known that $\frac{1}{n}\ln Z_n$ converges $\bb P$ a.s. towards $F(\beta,h)$.
The function $\eta\mapsto \ln Z_n$ is Lipschitz with constant $C=\beta\sqrt{n}$,
so if the disorder $\eta$ is gaussian, then one can use the concentration property of gaussian measures to deduce that
there exist constants $c_1$, $c_2$ such that
\begin{equation}\label{gaussian}
\forall x>0,\ \bb P(\frac{1}{n}\abs{\ln Z_n -\bb E\ln Z_n}>x)
\leq c_1 \exp(-c_2\frac{nx^2}{\beta^2}).
\end{equation}
This inequality remains valid if $\eta_0$ is a bounded variable or if the law of $\eta_0$ satisfies a log-Sobolev inequality (see \cite{Giacomin2007} for more details). We refer to \cite{Ledoux2001}
for a survey of the concentration measure phenomenon. Here we show that there is no need to make such assumptions to get exponential concentration of the free energy. 
Assuming only that $\bb E e^{\beta\abs{\eta}}<+\infty$, we prove that the inequality 
(\ref{gaussian}) is always valid for small values of $x$. 
%%%%%%%%
%thm concentration1
\begin{thm}\label{concentration1}
Assume that $\bb E e^{\beta\abs{\eta}}<+\infty$ and
set $K=2\max(e^{h+\lambda(\beta)},1)\times\max(e^{-h+\lambda(-\beta)},1)$. Then for all $n\geq 1$,
\begin{equation}\label{growth}
\bb E \exp(\pm t(\ln Z_n-\bb E\ln Z_n))\leq \exp(nKt^2)
\textrm{ for all }t\in [0,1]
\end{equation}
and \begin{equation}\label{speed}
\bb P(\pm\frac{1}{n}(\ln Z_n-\bb E \ln Z_n)>x)
\leq  \left\{\begin{aligned}
&\exp(-\ds\frac{nx^2}{4K}) &\textrm{ if } &x\in (0,2K],\\
&\exp(-n(x-K))  &\textrm{ if }&x\in
(2K, \infty).
\end{aligned}\right.
\end{equation}
\end{thm}
%%%%%%

We refer to \cite{CometsShigaYoshidaReview} for a review of directed polymers in
random environment.
In this model we consider $S=(S_n)_{n\in\bb N}$ the simple random walk on 
$\bb Z^d$ starting at $0$, defined on a probability space
$(\Sigma,\cal E,\textbf{P})$. Let $\eta=(\eta(n,x))_{(n,x)\in\bb N\times \bb
Z^d}$ be a sequence of real-valued, non-constant and i.i.d. random
variables defined on another probability space $(\Omega,\cal F, \bb P)$.
The path $S$ represents the directed polymer and $\eta$ the
random environment. For $n>0$, let $\Sigma_n$ be the set of paths of length $n$,
and define the random polymer measure
$\g P_n$ on the path space $(\Sigma,\cal E)$ by
\begin{equation}
\frac{d\g P_n}{d\textbf{P}}(S)=\ds\frac{1}{Z_n}\exp(\beta H_n(S)),
\end{equation}
where $\beta\in\bb R$ is the inverse temperature,
\begin{equation}
H_n(S)=\sum_{j=1}^n \eta(j,S_j),\textrm{ and } Z_n=\textbf{E}\exp(\beta H_n(S)).
\end{equation}
In other words, for a given realisation of the environment $\eta$, the measure $\g P_n$ gives to a polymer chain $S$ having an energy $-H_n(S)$ at temperature $T=\frac{1}{\beta}$, a weight proportional to $e^{\beta H_n(S)}$ (configurations of lowest energy are the most probable).
For simplicity we write $\bb E f(\eta)=\bb E f(\eta(0,0))$ for any $f$ such that $f\circ \eta(0,0)$ is integrable.
Let $\lambda(\beta)=\ln \bb E e^{\beta\eta}$ be the logarithmic moment generating function of
$\eta$. 
We shall use the point to point partition function:
\begin{equation}\label{part}
Z_{n}(x,y)=Z_{n}(x,y;\eta)=\textbf{E}^x \exp(\beta \sum_{j=1}^{n} \eta(j,S_j))\textbf{1}_{S_{n}=y},
\end{equation}
where $\textbf{E}^x$ is the expectation with respect to $\g P^x$, the measure of the random walk starting from $x$.
The Markov Property of the simple random walk yields that
\begin{equation}\label{markov}
Z_{n+m}(x,z)=\sum_{y}Z_{n}(x,y;\eta)Z_{m}(y,z;\Theta_n \eta),
\end{equation}
where $\Theta_n$ denotes the shift operator 
$$\Theta_{n}\eta(i,x)=\eta(i+n,x).$$
It is well known that this implies that the sequence $\bb E\ln Z_n$ is
superadditive, hence the limit
\begin{equation}\label{super}
p(\beta)=\Lim n \infty \frac{1}{n}\bb E\ln (Z_n)=
\sup_n \frac{1}{n}\bb E\ln (Z_n)\in (-\infty,\lambda(\beta)]
\end{equation}
exists. It is called the free energy of the polymer.
The function $\eta\mapsto \ln Z_n$ is Lipschitz with constant $C=\beta\sqrt{n}$,
so if the disorder $\eta$ is gaussian, then the concentration inequality (\ref{gaussian})
holds true (\cite{CarmonaHu2004}, Proposition 2.3), and therefore $\frac{1}{n}\ln Z_n$ converges $\bb P$ a.s. towards $p(\beta)$.
Using an inequality due to Lesigne and Volny (\cite{LesigneVolny}),
Comets, Shiga and Yoshida (\cite{CometsShigaYoshida}) proved that
if $\bb E e^{\beta\abs{\eta}}<+\infty$ for all $\beta>0$, then
for every $x>0$, there exists $n_0\in\bb N^*$ such that for any $n\geq n_0$,
\begin{equation}
\bb P(|\frac{1}{n}\ln Z_n-\frac{1}{n}\bb E\ln Z_n|>x)\leq
\exp(-\ds\frac{n^{\frac{1}{3}}x^{\frac{2}{3}}}{4}),
\end{equation}
which again allows them to prove that $\frac{1}{n}\ln Z_n$ converges $\bb P$ a.s. towards $p(\beta)$.
We improved this result in \cite{LiuWatbled} by showing that if $\bb E e^{\beta|\eta|}<\infty$
for a fixed $\beta>0$, then there exists $K>0$ such that for all $n\geq 1$,
\begin{equation*}%\label{speed}
\bb P(\pm\frac{1}{n}(\ln Z_n-\bb E \ln Z_n)>x)
\leq  \left\{\begin{aligned}
&\exp(-\ds\frac{nx^2}{4K}) &\textrm{ if } &x\in (0,2K],\\
&\exp(-n(x-K))  &\textrm{ if }&x\in
(2K, \infty).
\end{aligned}\right.
\end{equation*}
Here we extend this concentration property to functions of the form 
$\ln Y_n=\ln\textbf{E}f_n(S_n)\exp(\beta H_n(S))$, where $(f_n)$ is any sequence of bounded positive functions such that $\g P(f_n(S_n)>0)\neq 0$. 
%%%%%%%%%
%thm avecfn
\begin{thm}\label{avecfn}
Assume that $\bb E e^{\beta\abs{\eta}}<+\infty$ and
set $K=2\exp(\lambda(-\beta))+\lambda(\beta))$. Let $f_n$ be a sequence of
bounded positive functions on $\bb R^d$ such that for all $n\geq 1$, $\g P(f_n(S_n)>0)\neq 0$.
Set $Y_n=\g Ef_n(S_n)e^{\beta H_n(S)}$.
Then for all $n\geq 1$,
\begin{equation}\label{grfn}
\bb E e^{\pm t(\ln Y_n -\bb E\ln Y_n )}\leq 
\exp(nKt^2) \textrm{ for all }t\in [0,1],
\end{equation}
and
 \begin{equation}\label{spfn}%equation speedavecfn
\bb P(\pm\frac{1}{n}(\ln  Y_n-\bb E\ln Y_n)>x)
\leq  \left\{\begin{aligned}
&\exp(-\ds\frac{nx^2}{4K}) &\textrm{ if } &x\in (0,2K],\\
&\exp(-n(x-K))  &\textrm{ if }&x\in
(2K, \infty).
\end{aligned}\right.
\end{equation}
\end{thm}
%%%%%%

In \cite{CarmonaHu2004}, Carmona and Hu consider a gaussian environment and show among other things a large deviations result for the end position $S_n$ of the polymer under the polymer measure
$\g P_n$ (\cite{CarmonaHu2004}, Theorems 1.1 and 1.2). Their proof relies on the concentration inequality (\ref{gaussian}). Here we prove that their large deviations result holds for every environment such that $\bb E e^{\beta|\eta|}<\infty$.
Let $B_d=\{x\in\bb R^d; \N{x}_1=\abs{x_1}+\cdots+\abs{x_d}\leq 1\}$ be the closed unit ball of $\bb R^d$ in the $l^1$-norm, $\mathring{B_d}$ be the corresponding open ball.
%%%%%%%%%%%%%%%%
%thm deviation
\begin{thm}\label{deviation}
Let $\beta >0$ and assume that $\bb E e^{\beta\abs{\eta}}<+\infty$. 
Then there exists a convex rate function $I_\beta:B_d\to [0,\ln(2d)+p(\beta)]$ such that $\bb P$ a.e.,
\begin{align}
\limsup_{n\to \infty}\frac{1}{n} \ln\mathbf{P}_n (\frac{S_n}{n}\in F)&\leq -\inf_{x\in F} I_\beta(x) \textrm{ for } F \textrm{ closed }\subset B_d, \label{ratesup}\\ %equation ratesup
\liminf_{n\to \infty}\frac{1}{n} \ln\g P_n(\frac{S_n}{n}\in G)&\geq -\inf_{x\in G} I_\beta(x) \textrm{ for } G \textrm{ open }\subset B_d. \label{rateinf}%equation rateinf
\end{align}
\end{thm}
%%%%%%%%
%thm pointwise
\begin{thm}\label{pointwise}
Let $\beta >0$ and assume that $\bb E e^{\beta\abs{\eta}}<+\infty$. 
Then for every $x\in \mathring{B_d}\cap \bb Q^d$,
\begin{equation}\label{point}%equation point
\lim_{n\to \infty}\frac{1}{n} \ln\g P_n (\frac{S_n}{n}=x)= - I_\beta(x)\   \bb P \ a.s.,
\end{equation}
where we take the limit along $n$ such that $\g P(S_n=nx)>0$. 
Moreover, $I_\beta(0)=0$,
$I_\beta(x_1,\cdots,x_d)=I_\beta(\pm x_{\sigma(1)},\cdots,\pm x_{\sigma(d)})$ for any permutation $\sigma$ of $\{1,\cdots,d\}$, and for $e_1=(1,0,\cdots, 0)\in \bb Z^d$, we have
\begin{equation}\label{e1}%equation e1
I_\beta(e_1)=\ln (2d)+p(\beta).
\end{equation}
\end{thm}
%%%%%%%%%

The paper is organized as follows.
In a first part we recall in a simpler form a theorem of \cite{LiuWatbled}, which gives optimal conditions 
on the martingale $M_n$ to provide exponential bounds for $\bb P(\frac{M_n}{n} > x)$.
 In the second part we prove Theorem \ref{concentration1} by writing $\ln Z_n-\bb E \ln Z_n$ as a
sum of $(\cal F_j)_{1\leq j\leq n}$ martingale differences, 
where $\cal F_j=\sigma(\eta_i: 1\leq i\leq j)$, and showing that the conditional exponential moments of the martingale differences are uniformly bounded.
In the third part we prove in a similar manner Theorem \ref{avecfn} by writing $\ln Y_n-\bb E \ln Y_n$ as a sum of $(\cal F_j)_{1\leq j\leq n}$ martingale differences, using this time the filtration 
$\cal F_j=\sigma(\eta(i,x): 1\leq i\leq j, x\in \bb Z^d)$.
In the last part we explain how to use our concentration results (\ref{spfn}) of Theorem \ref{avecfn}
to prove Theorem \ref{deviation} and Theorem \ref{pointwise}.

%%%%%%% section 
\section{Exponential inequalities for martingales}
\setcounter{equation}{0}

Let $(\Omega,\cal F,\bb P)$ be a probability
space and let $\cal F_0=\{\emptyset,\Omega\}\subset \cal
F_1\subset\cdots\subset \cal F_n$ be an increasing sequence of
sub-$\sigma$-fields of $\cal F$. Let $ X_1, ..., X_n$ be a sequence
of real-valued martingale differences defined on $(\Omega,\cal
F,\bb P)$, adapted to the filtration $(\cal F_k)$: that is, for each
$1\leq k \leq n$,   $X_k$ is $\cal F_k$ measurable and $\bb
E(X_k|\cal F_{k-1}) = 0$. Set
\begin{equation}
 M_n = X_1 + ... + X_n.
 \end{equation}
 If the martingale differences $X_i$ are uniformly bounded by a constant $a$,
 then Azuma-Hoeffding's inequality (\cite{Hoeffding},\cite{Azuma})
 states that
 \begin{equation}\label{AH}
\bb P(\frac{M_n}{n} > x) \leq   e^{- \frac{nx^ 2}{a^ 2}}.
\end{equation}
Lesigne and Voln\'y
proved (\cite{LesigneVolny}) that if for some constant $K>0$ and all $k=1, ..., n$,
 \begin{equation}
 \label{exp-moment}
       \bb E e^{ |X_k| } \leq K,
\end{equation}
then for any $x>0$,
\begin{equation}
\label{LV}
 \bb P(\frac{M_n}{n} > x) = O ( e^{ - \frac{1}{4}  x^{2/3}
 n^{1/3}}).
\end{equation}
They also showed  that this is the best possible inequality that we
can have under the condition (\ref{exp-moment}),  even in the class
of stationary and ergodic sequences of martingale differences, in
the sense that there exist such sequences of martingale differences
$(X_i)$ satisfying (\ref{exp-moment}) for some $K>0$, but
\begin{equation}
\label{LV2}
\bb P(\frac{M_n}{n} > 1)> e^{ - c n^{1/3}}
 \end{equation}
for some constant $c>0$ and infinitely many $n$.
We showed in \cite{LiuWatbled} that the best condition for having an exponential inequality of the form
$$\bb P(\frac{M_n}{n} > x) \leq Ce^{-c(x)n}$$
is to replace the
expectation in (\ref{exp-moment}) by the conditional one given $\cal
F_{k-1}$.
We recall rapidly one of our results, presented in a slightly simpler form, and its proof
(\cite{LiuWatbled}, Theorem 2.1).
%%%%%%%%%
\begin{thm}\label{superdiff}
Let $(X_i)_{1\leq i\leq n}$ be a finite sequence of martingale
differences. If for some constant $K>0$ and all $i=1, ..., n$,
\begin{equation}\label{K}
\bb E(e^{\abs{X_i}}|\cal F_{i-1})\leq K\ \ \ a.s.,
\end{equation}
then:
\begin{equation}\label{growthsuperdiff}
 \bb E e^{\pm tM_n}\leq \exp(nKt^2)
\textrm{ for all }t\in [0,1],
\end{equation}
and 
\begin{equation}\label{concsuperdiff}
\bb P(\frac{\pm M_n}{n}>x)\leq \left\{\begin{aligned}
&\exp(-\ds\frac{nx^2}{4K}) &\textrm{ if } &x\in (0,2K],\\
&\exp(-n(x-K))  &\textrm{ if }&x\in
(2K, \infty).
\end{aligned}\right.
\end{equation}
\end{thm}
%%%%%%

We need the following simple lemma.
%%%%%%
%lemsuperdiff
\begin{lem}\label{lemsuperdiff} Let $X$ be a real-valued random variable defined
on some probability space $(\Omega, \cal F, \bb P)$, with $\bb E X\leq 0$ and
$ \bb E e^{\abs{X}}\leq K$ for some $K>0$.
Then for all $t\in [0,1]$,
\begin{equation}
\bb E e^{tX}\leq \exp(Kt^2).
\end{equation}
\end{lem}
%%%%%%
%preuve lemuperdiff
\textbf{Proof.} Let $t\in [0,1]$. Since $\bb E X\leq 0$, we have:
$$\begin{aligned}
\bb E e^{tX} 
&=\sum_{k=0}^\infty t^k\bb E\frac{X^k}{k!}
\leq 1+ \sum_{k=2}^\infty t^k\bb E\frac{X^k}{k!} \\
&\leq 1+t^2\sum_{k=2}^\infty \bb E\frac{\abs{X}^k}{k!}
\leq 1+t^2\bb E e^{\abs{X}}
\leq 1+Kt^2\leq \exp(Kt^2) .
\end{aligned}$$
   \hfill{\rule{2mm}{2mm}\vskip3mm \par}

%preuve thmsuperdiff
\textbf{Proof of Theorem \ref{superdiff}.}
By Lemma \ref{lemsuperdiff}, we obtain that for every $i$ and for every $t\in (0,1)$,
$$\bb E (e^{t X_i}|\cal F_{i-1})\leq \exp(Kt^2) \  a.s.$$
By induction on $n$ we get immediately (\ref{growthsuperdiff}) for $+M_n$.
Then $\forall x>0$, $\forall t\in [0,1]$,
$$\bb P(\frac{M_n}{n}>x)=\bb P(e^{tM_n}> e^{tnx})
\leq e^{-ntx}\bb E e^{tM_n}\leq \exp(-n(tx-Kt^2)).$$
As 
\begin{equation}\label{optim}
 \sup_{t\in [0,1]}(tx-Kt^2)= \left\{\begin{aligned}
&\frac{x^2}{4K} &\textrm{ if } &x\in (0,2K],\\
&x-K  &\textrm{ if }&x\in
(2K, \infty),
\end{aligned}\right.
\end{equation}
we obtain (\ref{concsuperdiff}) for $+M_n$. 
We apply the same argument to the sequence $(-X_i)$ and obtain the full inequalities
(\ref{growthsuperdiff}) and (\ref{concsuperdiff}).\hfill{\rule{2mm}{2mm}\vskip3mm \par}

%%%%%%% Sect 
\section{Disordered pinning model: proof of Theorem \ref{concentration1}}
\setcounter{equation}{0}

We write $\ln Z_n-\bb E\ln Z_n$ as a
sum of $(\cal F_j)_{1\leq j\leq n}$ martingale differences:
\begin{equation}
\ln Z_n-\bb E\ln Z_n=\ds\sum_{j=1}^nV_{n,j}, \quad \textrm{ with }
V_{n,j}=\bb E_j\ln Z_n-\bb E_{j-1}\ln Z_n,
 \end{equation}
 where $\bb E_j$
denotes the conditional expectation with respect to $\bb P$ given
$\cal F_j$, $\cal F_j=\sigma(\eta_i: 1\leq i\leq j)$.
Theorem \ref{concentration1} is then a direct consequence of Theorem \ref{superdiff}.
All we have to do is to check condition (\ref{K}), which is done in the following lemma.

%lemme CSY
\begin{lem}\label{CSY} For every $1\leq j\leq n$, we have:
\begin{equation}\label{CSYbis}
\bb E_{j-1}\exp(\abs{V_{n,j}})\leq
K:= 2\max(e^{h+\lambda(\beta)},1)\times\max(e^{-h+\lambda(-\beta)},1).
\end{equation}
\end{lem}

%preuve lemme CSY
\textbf{Proof.}
Let us define
 \begin{equation}
 Z_{n,j}=\textbf{E}\exp(\sum_{k=1,k\neq j}^n(\beta \eta_k+h)\delta_k)\delta_n.
  \end{equation}
Since $\bb E_{j-1}\ln Z_{n,j}=\bb E_j\ln Z_{n,j},$ we have:
\begin{equation}\label{W}
V_{n,j}=\bb E_j\ln \frac{Z_n}{Z_{n,j}}-
\bb E_{j-1}\ln \frac{Z_n}{Z_{n,j}},
\end{equation}
therefore
\begin{equation}
\bb E_{j-1}\exp(V_{n,j})=
\exp(-\bb E_{j-1} \ln \frac{Z_n}{Z_{n,j}})
\bb E_{j-1}\exp(\bb E_j\ln\frac{Z_n}{Z_{n,j}}).
\end{equation}
Using Jensen's inequality and the fact that $\cal F_{j-1}\subset\cal F_j$, we get:
\begin{equation}\label{CSY1}
\bb E_{j-1}\exp(V_{n,j})\leq
\bb E_{j-1}(\frac{Z_n}{Z_{n,j}})^{-1}
\bb E_{j-1}\frac{Z_n}{Z_{n,j}}.
\end{equation}
Now let us write
\begin{equation}
\frac{Z_n}{Z_{n,j}}=\alpha_0+\alpha_1\exp(\beta\eta_j+h),
\end{equation}
with \begin{equation}
\alpha_l=\frac{\g E\exp(\sum_{k=1,k\neq j}^n(\beta \eta_k+h)\delta_k)\delta_n  \textbf{1}_{\delta_j=l}}{Z_{n,j}}\ \textrm{ for }l=0,1.
\end{equation} 
We consider the $\sigma$-algebra $\cal F_{n,j}=\sigma(\eta_k;1\leq k\leq n, k\neq j)$. 
Then $\cal F_{j-1}\subset\cal F_{n,j}$ and, therefore,
\begin{align*}
\bb E_{j-1}\frac{Z_n}{Z_{n,j}}
&=\bb E_{j-1}\bb E(\alpha_0+\alpha_1\exp(\beta\eta_j+h)\mid \cal F_{n,j})\\
&=\bb E_{j-1}(\alpha_0+\alpha_1\exp(\lambda(\beta)+h))\\
&=\bb E_{j-1}\alpha_0+\bb E_{j-1}(\alpha_1)\exp(\lambda(\beta)+h).
\end{align*}
As $\alpha_0+\alpha_1=1$ we deduce that
\begin{equation}\label{prem}
\bb E_{j-1}\frac{Z_n}{Z_{n,j}}\leq \max(e^{h+\lambda(\beta)},1).
\end{equation}
For the same reason, and using the convexity of the inverse function, we obtain
\begin{align*}
( \frac{Z_n}{Z_{n,j}})^{-1}&=(\alpha_0+\alpha_1\exp(\beta\eta_j+h))^{-1}\\
&\leq \alpha_0+\alpha_1\exp(-\beta\eta_j-h).
\end{align*}
The same reasoning as before then leads to
\begin{equation}\label{deux}
\bb E_{j-1}( \frac{Z_n}{Z_{n,j}}) ^{-1}\leq \max(e^{-h+\lambda(-\beta)},1).
\end{equation}
Combining (\ref{CSY1}), (\ref{prem}) and (\ref{deux}) we get the inequality
\begin{equation}
\bb E_{j-1}\exp(V_{n,j})\leq \max(e^{h+\lambda(\beta)},1)\max(e^{-h+\lambda(-\beta)},1).
\end{equation}
The same inequality holds with $-V_{n,j}$ instead of $V_{n,j}$, from which we deduce (\ref{CSYbis}).
\hfill{\rule{2mm}{2mm}\vskip3mm \par}

We conclude this section by noticing that the concentration inequality (\ref{speed})
implies immediately the following convergence result. 
%%%%%%%%
%cor convLp
\begin{cor}\label{convLp}
Assume that $\bb E e^{\beta\abs{\eta}}<+\infty$.
Then:
\begin{equation}\label{convergenceenergy}
\frac{\ln Z_n}{n}-\bb E\frac{\ln Z_n}{n}\to 0 \ \ a.s. \textrm{ and in } L^p, p\geq 1.
\end{equation}
\end{cor}

%%%%%%% Sect 
\section{Directed polymers in random environment: proof of Theorem \ref{avecfn}}
\setcounter{equation}{0}

The proof of Theorem \ref{avecfn} follows the same line as the proof of Theorem \ref{concentration1}.
We write $\ln Y_n-\bb E\ln Y_n$ as a
sum of $(\cal F_j)_{1\leq j\leq n}$ martingale differences:
\begin{equation}
\ln Y_n-\bb E\ln Y_n=\ds\sum_{j=1}^nV_{n,j}, \quad \textrm{ with }
V_{n,j}=\bb E_j\ln Y_n-\bb E_{j-1}\ln Y_n,
 \end{equation}
 where this time $\bb E_j$
denotes the conditional expectation with respect to $\bb P$ given
$\cal F_j$, $\cal F_j=\sigma(\eta(i,x): 1\leq i\leq j, x\in \bb Z^d)$.
According to Theorem \ref{superdiff}, to prove Theorem \ref{avecfn} we only have to bound the conditional exponential moments of the martingale differences,
which is done in the following lemma.

%lemme CSYfn
\begin{lem}\label{CSYfn} For every $1\leq j\leq n$, we have:
\begin{equation}\label{inegb}
\bb E_{j-1}\exp(\abs{V_{n,j}})\leq
K:= 2\exp(\lambda(\beta)+\lambda(-\beta)).
\end{equation}
\end{lem}

%preuve lemme CSYfn
\textbf{Proof.} For every $1\leq j\leq n$, we define
 \begin{equation}
  H_{n,j}(S)=\sum_{1\leq k\leq n,k\neq j}\eta(k,S_k) \textrm{ and } 
 Y_{n,j}=\textbf{E}f_n(S_n)\exp(\beta H_{n,j}(S)).
  \end{equation}
Since $\bb E_{j-1}\ln Y_{n,j}=\bb E_j\ln Y_{n,j},$ we obtain as in Lemma \ref{CSY}
that
\begin{equation}\label{CSYfn1}
\bb E_{j-1}\exp(V_{n,j})\leq
\bb E_{j-1}(\frac{Y_n}{Y_{n,j}})^{-1}
\bb E_{j-1}\frac{Y_n}{Y_{n,j}}.
\end{equation}
Now let us write
\begin{equation}
\frac{Y_n}{Y_{n,j}}=\sum_{x\in\bb Z^d}\alpha_x\exp(\beta\eta(j,x)),
\end{equation}
with 
 \begin{equation}\label{def}
\alpha_x=\frac{\g E f_n(S_n)e^{\beta H_{n,j}(S)} \textbf{1}_{S_j=x}}{Y_{n,j}}.
  \end{equation}
We consider the $\sigma$-algebra $\cal F_{n,j}=\sigma(\eta(k,x);1\leq k\leq n, k\neq j, x\in \bb Z^d)$. 
Using that $\cal F_{j-1}\subset\cal F_{n,j}$, the $\alpha_x$ are $\cal E_{n,j}$-measurable, and the
$\exp(\beta\eta(j,x))$ are independent of $\cal E_{n,j}$, we obtain:
\begin{align*}
\bb E_{j-1}\frac{Y_n}{Y_{n,j}}
&=\bb E_{j-1}\bb E(\sum_{x\in\bb Z^d}\alpha_x\exp(\beta\eta(j,x))\mid \cal F_{n,j})\\
&=\bb E_{j-1}\sum_{x\in\bb Z^d}\alpha_x\exp(\lambda(\beta)).
\end{align*}
As $\sum_{x\in\bb Z^d}\alpha_x=1$, we get
\begin{equation}\label{CSYfn2}
\bb E_{j-1}\frac{Y_n}{Y_{n,j}}=\exp(\lambda(\beta)).
\end{equation}
The inverse function is convex and $\sum_{x\in\bb Z^d}\alpha_x=1$, therefore
\begin{align*}
( \frac{Y_n}{Y_{n,j}})^{-1}&=(\sum_{x\in\bb Z^d}\alpha_x\exp(\beta\eta(j,x)))^{-1}\\
&\leq\sum_{x\in\bb Z^d}\alpha_x\exp(-\beta\eta(j,x)).
\end{align*}
The same reasoning as before then leads to
\begin{equation}\label{CSYfn3}
\bb E_{j-1}( \frac{Y_n}{Y_{n,j}}) ^{-1}\leq \exp(\lambda(-\beta)).
\end{equation}
Combining (\ref{CSYfn1}), (\ref{CSYfn2}) and (\ref{CSYfn3}) we get the inequality
\begin{equation}
\bb E_{j-1}\exp(V_{n,j})\leq \exp(\lambda(\beta)+\lambda(-\beta)).
\end{equation}
The same inequality holds with $-V_{n,j}$ instead of $V_{n,j}$, from which we deduce (\ref{inegb}).
\hfill{\rule{2mm}{2mm}\vskip3mm \par}

As in the preceding section, we notice that the concentration inequality (\ref{spfn})
implies immediately that $\frac{\ln Y_n}{n}-\bb E\frac{\ln Y_n}{n}$ converges to $0$ a.s. and also in $L^p$ for all $p\geq 1$. 
In particular, when $f_n=1$, we conclude that $\frac{\ln Z_n}{n}$ converges
$\bb P$ a.s. towards $p(\beta)$.
In the same manner, if the sequence $\bb E\ln Y_n$ is superadditve, then it converges towards a limit $L(\beta)$, and the concentration inequality (\ref{spfn})
implies that $\frac{\ln Y_n}{n}$ converges
$\bb P$ a.s. towards $L(\beta)$.
Let us denote by $\g E_n f=\int fd\g P_n$ the expectation of a function $f$ with respect to the polymer measure $\g P_n$.
In particular, with the notations of Theorem \ref{avecfn}, we have $\g E_n f_n(S_n)=\frac{Y_n}{Z_n}$,
and we can deduce from Theorem \ref{avecfn} the following result,
which will be used repeatedly in the next section.
 
%%%%%%
%cor corfn
\begin{cor}\label{corfn}
Assume that $\bb E e^{\beta\abs{\eta}}<+\infty$. Let $(f_n)$ be a sequence of
bounded positive functions on $\bb R^d$ such that for all $n\geq 1$, $\g P(f_n(S_n)>0)\neq 0$.
Set $Y_n=\g E f_n(S_n)e^{\beta H_n(S)}$.
Then:
\begin{equation}\label{convergenceenergyfn}
\frac{\ln Y_n}{n}-\bb E\frac{\ln Y_n}{n}\to 0 \ \ a.s. \textrm{ and in } L^p \textrm{ for all }p\geq 1.
\end{equation}
If moreover the sequence $\bb E\ln Y_n$ is superadditive, then the limit
\begin{equation}\label{limitQ}
L(\beta)=\limf{n}{+\infty}\frac{1}{n}\bb E\ln Y_n=\sup_{n\geq 1}\frac{1}{n}\bb E\ln Y_n
\end{equation}
exists, so does the limit 
\begin{equation}
\limf{n}{+\infty}\frac{1}{n}\bb E\ln \g E_n f_n(S_n)=L(\beta)-p(\beta),
\end{equation}
and then:
\begin{equation}
\limf{n}{+\infty}\frac{1}{n}\ln \g E_n f_n(S_n)=L(\beta)-p(\beta)\ \ a.s. \textrm{ and in } L^p, \textrm{ for all }p\geq 1.
\end{equation}
\end{cor}

%%%%%%% Sect 
\section{Large deviations for directed polymers in random environment}
\setcounter{equation}{0}

Theorems \ref{deviation} and \ref{pointwise} are proved in
\cite{CarmonaHu2004} (Theorems 1.1 and 1.2) in the case where 
$(\eta(n,x))_{(n,x)\in\bb N\times\bb Z^d}$ are i.i.d. $\cal N(0,1)$ gaussian random variables.
Their proofs rely essentially on subadditivity, and
on the concentration property of the gaussian measure.
We can extend it to the case of a general environment, assuming only that 
$\bb E e^{\beta\abs{\eta}}<+\infty$, because the subaddivity arguments are still valid,
and Theorem \ref{avecfn} provides the concentration properties we need.
More precisely
we shall use several times (that is for several different sequences $(f_n)$'s)
Corollary \ref{corfn}.
Without loss of generality we shall assume that $\bb E\eta=0$.
Our proof of Theorems \ref{deviation} and \ref{pointwise} are essentially the same as
Theorems 1.1 and 1.2 of \cite{CarmonaHu2004}, so we won't detail them. Instead let us 
give two lemmas which illustrate how Theorem \ref{avecfn} is used.
In Lemma 3.1 Carmona and Hu establish that for any $\lambda >0$, for any $x\in B_d$ and any sequence $x_n\in\bb R^d$ satisfying $x_n/n\to x$, the following limit exists thanks to subadditivity:
\begin{equation}
\underset{n\to +\infty}{\lim}\frac{\bb E\ln \g E e^{-\lambda\N{S_n-x_n}_1 }e^{\beta H_n(S)}}{n}=\phi_\lambda(x).
\end{equation}
We first consider $f_n(z)=e^{-\lambda\N{z-x_n}_1 }$, where $\lambda>0$ and $x_n\in \bb R^d$,
and applying Corollary \ref{corfn} we obtain:

%lemma lemexp
\begin{lem} \label{lemexp} For any $x\in B_d$ and any sequence $x_n\in\bb R^d$ satisfying $x_n/n\to x$, the following limits exist a.s. 
and in $L^p$, for all $p\geq 1$:
\begin{equation}\label{cexp}
\underset{n\to +\infty}{\lim}-\frac{1}{n}\ln \g E_n e^{-\lambda\N{S_n-x_n}_1 }=p(\beta)-\phi_\lambda(x).
\end{equation}
\end{lem}

The function $I_\beta^{(\lambda)}=p(\beta)-\phi_\lambda(x)$ is nondecreasing in $\lambda$. The function 
$$I_\beta=\sup_{\lambda>0}I_\beta^{(\lambda)}$$
is then convex and lower semicontinous on $B_d$ with values in $[0,p(\beta)+\ln(2d)]$,
and this is the rate function which existence is claimed in Theorem \ref{deviation}.
We consider now successively $f_n(z)=\mathbf{1}_{\{\N{z-nx}_{1}\leq n\eps\}}$, where $x\in B_d$ and $\eps >0$, and $f_n(z)=\mathbf{1}_{\{\N{z-nx}_{1}< n\eps\}}$.

%lemma lemboule
\begin{lem} \label{lemboule}
For any $x\in B_d$ and any $\eps>0$, the following limits exist:
\begin{equation}\label{E2}
L_\beta(x,\eps)=\limf{n}{+\infty} \frac{1}{n}\bb E\ln \g E_n \mathbf{1}_{\{\N{S_n-nx}_{1}\leq n\eps\}},
\end{equation}
and
\begin{equation}\label{cbou}
L_\beta(x,\eps)=\limf{n}{+\infty} \frac{1}{n}\ln \g E_n \mathbf{1}_{\{\N{S_n-nx}_{1}\leq n\eps\}} \ \ a.s. 
\textrm{ and in } L^p, \textrm{ for all } p\geq 1,
\end{equation}
as well as
\begin{equation}\label{E12}
\mathring{L}_\beta(x,\eps)=\limf{n}{+\infty} \frac{1}{n}\bb E\ln \g E_n\mathbf{1}_{\{\N{S_n-nx}_{1}< n\eps\}},
\end{equation}
and
\begin{equation}\label{E13}
\mathring{L}_\beta(x,\eps)=\limf{n}{+\infty} \frac{1}{n}\ln \g E_n\mathbf{1}_{\{\N{S_n-nx}_{1}< n\eps\}} \ \ a.s. 
\textrm{ and in } L^p, \textrm{ for all } p\geq 1.
\end{equation}
\end{lem}
%preuve lem lemboule
\textbf{Proof of Lemma \ref{lemboule}.}
Let us fix $x\in B_d$ and $\eps >0$. 
We set $$Y_n=\g E\mathbf{1}_{\{\N{S_n-nx}_{1}\leq n\eps\}}e^{\beta H_n(S)},
\ v_n=\bb E\ln Y_n,$$
and show that $v_n$ is superadditive. The method of proof is the same as for proving that 
$\bb E\ln Z_n$ is superadditive.
If $\N{S_n-nx}_{1}\leq n\eps$ and $\N{(S_{n+m}-S_n)-mx}_{1}\leq m\eps$
then $\N{S_{n+m}-(n+m)x}_{1}\leq (n+m)\eps$,
so we have
\begin{equation}
\mathbf{1}_{\{\N{S_{n+m}-(n+m)x}_{1}\leq (n+m)\eps\}}\geq
\mathbf{1}_{\{\N{S_{n}-nx}_{1}\leq n\eps\}}\times
\mathbf{1}_{\{\N{S_{n+m}-S_n -  mx}_{1}\leq m\eps\}}.
\end{equation}
Therefore
\begin{equation}\label{boule1}
Y_{n+m}\geq
 \g E\mathbf{1}_{\{\N{S_n-nx}_{1}\leq n\eps\}}e^{\beta H_n(S)}\times
\mathbf{1}_{\{\N{S_{n+m}-S_n -  mx}_{1}\leq m\eps\}}
e^{\beta\sum_{i=n+1}^{n+m} \eta(i,S_i)}.
\end{equation}
The right handside is equal to
\begin{equation}\label{boule2}
\sum_{y\in\bb Z^d} \g E\mathbf{1}_{\{\N{S_n-nx}_{1}\leq n\eps\}}e^{\beta H_n(S)}\mathbf{1}_{\{S_n=y\}}
\g E^y\mathbf{1}_{\{\N{S_{m}-mx}_{1}\leq m\eps\}}e^{\beta\sum_{i=1}^{m} \eta(n+i,S_i)}.
\end{equation}
Let $\sigma_n$  be the probability measure defined on $\Sigma$ by 
\begin{equation}
\frac{d\sigma_n}{d\g P}(S)=\frac{\mathbf{1}_{\{\N{S_n-nx}_{1}\leq n\eps\}}e^{\beta H_n(S)}}
{Y_n}.
\end{equation}
Then (\ref{boule1}) and (\ref{boule2}) imply
\begin{equation}
Y_{n+m}\geq
\sum_{y\in\bb Z^d} \sigma_n(S_n=y)Y_n
\g E^y\mathbf{1}_{\{\N{S_{m}-mx}_{1}\leq m\eps\}}
e^{\beta\sum_{i=1}^{m} \eta(n+i,S_i)}.
\end{equation}
Using the concavity of the logarithm we obtain 
\begin{equation}
\ln Y_{n+m}\geq \ln Y_n +
\sum_{y\in\bb Z^d} \sigma_n(S_n=y)
\ln \g E^y\mathbf{1}_{\{\N{S_{m}-mx}_{1}\leq m\eps\}}
e^{\beta\sum_{i=1}^{m}\eta(n+i,S_i)}.
\end{equation}
We take the conditional expectation with respect to $\bb P$ given 
$\cal F_n$ (recall that $\bb E_n$
denotes the conditional expectation with respect to $\bb P$ given
$\cal F_n$, $\cal F_n=\sigma(\eta(i,x): 1\leq i\leq n, x\in \bb
Z^d$), and we obtain the following inequality
\begin{equation}
\bb E_n\ln Y_{n+m}\geq
\ln Y_n +\bb E\ln Y_m.
\end{equation}
Integrating with respect to $\bb P$ we conclude that $(v_n)$ is superadditive:
$$v_{n+m}\geq v_n+v_m.$$
Therefore the limit
\begin{equation}
l_\beta(x,\eps)=\limf{n}{+\infty}\frac{v_n}{n}=\sup_{n\geq 1}\frac{v_n}{n}
\end{equation}
exists, and so does, by Corollary \ref{corfn}, the limit
\begin{equation}
L_\beta(x,\eps)=\limf{n}{+\infty} \frac{1}{n}\ln \g E_n \mathbf{1}_{\{\N{S_n-nx}_{1}\leq n\eps\}}
\ \bb P\  a.s. \textrm{ and in } L^p \textrm{ for all }p\geq 1,
\end{equation}
with
$$L_\beta(x,\eps)=l_\beta(x,\eps)-p(\beta).$$
The proof of \eqref{E12} and \eqref{E13} is the same as for \eqref{E2} and \eqref{cbou}.
\hfill{\rule{2mm}{2mm}\vskip3mm \par}

\medskip

\textbf{Acknowledgements.}
I am very grateful to an anonymous referee who made a lot of 
very valuable comments. I would like to thank also Quansheng Liu for fruitful discussions.

\bibliographystyle{plain}
\bibliography{biblioprobconc}

\end{document}